\documentclass[11pt,reqno]{amsart}
\usepackage{amssymb,latexsym}
\usepackage{graphicx}
\usepackage{enumerate,multicol}
\usepackage{fancyhdr,amssymb}
\usepackage{url}        

\theoremstyle{plain} \numberwithin{equation}{section}
\newtheorem{thm}{Theorem}[section]
\newtheorem{theorem}[thm]{Theorem}

\newtheorem{corollary}[thm]{Corollary}

\newtheorem{proposition}[thm]{Proposition}

\begin{document}

\fancyhead{}
\renewcommand{\headrulewidth}{0pt}

\setcounter{page}{1}

\title[On The Horadam Symbol Elements]{On The Horadam Symbol Elements}
\author{Sai Gopal Rayaguru}
\address{ National Institute of Technology\\
                Rourkela, India}
\email{saigopalrs@gmail.com}
\author{Diana Savin}
\address{ Ovidius University,Bd. Mamaia 124, 900527,\\
                Constanta, Romania}
\email{savin.diana@univ-ovidius.ro; dianet72@yahoo.com }
\author{Gopal Krishna Panda}
\address{ National Institute of Technology\\
                Rourkela, India}
\email{gkpanda\_nit@rediffmail.com}

\begin{abstract}
 Horadam symbol elemnts are introduced. Certain properties of these elements are explored. Some well known identities such as Catalan identity, Cassini formula and d'Ocagne's identity are obtained for these elements. In the last section we use these properties for to find zero divisors in symbol algebras over cyclotomic fields of finite fields.
\end{abstract}

\maketitle \textbf{Key words:} Recurrence Relations, Quaternions, Octonions, Symbol algebra.   

\maketitle \textbf{2010 Subject classification [A.M.S.]:} 11R52;11B37; 11B39; 20G20.

\bigskip

\section{Introduction}

\bigskip

\bigskip

Let $N$ be a positive integer and $K,$ a field such that
$char(K)$ does not divide $N$ and it contains a primitive $N^{th}$ root of unity: $\omega$. Let
$K^* = K \setminus\{0\},~a, b \in K^*$, and let $S$ be the algebra over K generated by
elements $x$ and $y$, with
$x^N = a$, $y^N = b$, $yx = \omega xy.$
This algebra is called a symbol algebra (also known as a power norm
residue algebra) and is denoted by  $\big(\frac{a,b}{K,\omega}\big).$
Observe that $S$ has a $K$-basis $\{x^i y^j :0 \leq i, j < N\}.$ The symbol algebra of degree 2 is called the quaternion algebra and here the symbol algebras are generalization of the quaternion algebras.\\ 

In 1963, Horadam\cite{Horadam} defined the $n$-th Fibonacci quaternion $Q_n$ as $$Q_n=F_n+iF_{n+1}+jF_{n+2}+kF_{n+3},$$ where $F_n$ denote the $n$-th Fibonacci number and $i,j,k$ are the quaternion units, satisfying $$i^2=j^2=k^2=-1,ij = -ji = k, jk = -kj = i, ki = -ik =j.$$ Subsequently many authors studied the extension and generalization of the Fibonacci quaternion and investigated their properties. Flaut and Shpakivskyi\cite{Flaut2013a} studied on the properties for the generalized Fibonacci quaternions and Fibonacci-Narayana quaternions in a generalized quaternion algebra. Variants of octonion algebra involving generalized Fibonacci octonions were also studied by many authors (see \cite{Gamaliel, Halici2012,Halici2013, Halici2016a, Ipek, Halici2016, Kosal2017, Savin, Savin2016, Yilmaz2017}).\\
 
 Flaut and Savin\cite{Flaut2012} and after Flaut, Savin and Iorgulescu \cite{Flaut2013b} considered the symbol algebra of degree 3 and defined the $n$-th Fibonacci symbol elements $\mathit{F}_{n}$ and $n$-th Lucas symbol elements $\mathit{L}_{n}$ as\\
  
$\mathit{F}_{n}=f_{n}.1+f_{n+1}.x+f_{n+2}.x^2+f_{n+3}.y+f_{n+4}.xy+f_{n+5}.x^2y+f_{n+6}.y^2+f_{n+7}.xy^2+f_{n+8}.x^2y^2$\\

 and \\
 
$\mathit{L}_{n}=l_{n}.1+l_{n+1}.x+l_{n+2}.x^2+l_{n+3}.y+l_{n+4}.xy+l_{n+5}.x^2y+l_{n+6}.y^2+l_{n+7}.xy^2+l_{n+8}.x^2y^2$\\ respectively, where $f_n$ and $l_n$ are the $n$-th Fibonacci and and $n^{th}$ Lucas number.\\  

Although, the notation used for defining the Fibonacci symbol elements is same as the notation for the Fibonacci numbers in the Fibonacci quaternion, we confirm their distinct behavior. We believe that the readers will not be confused with these notations. Flaut and Savin\cite{Flaut2017} studied some properties and applications of $(a,b,x_0,x_1)$-numbers and quaternion. Savin\cite{Savin2017} defined and studied properties of special numbers, special quaternions and special symbol elements. In a recent paper, Flaut and Savin\cite{Flaut2018} presented some applications of special numbers obtained from a difference equation of degree three. We generalize the concept of Fibonacci and Lucas symbol elements to the Horadam symbol elements. Moreover, we study the  Binet formula, generating function and some identities involving Horadam symbol elements. We will show that these elements satisfy many properties similar to that of Horadam numbers.\\

The Horadam sequence $\{w_{k}(a_0,a_1,p,q) \},{k \geq 1}$ is defined by 
\begin{equation}{\label{1.1}}
\hspace{3.0cm} w_{k+1} = pw_{k} + qw_{k-1};w_{0} = a_0,w_{1} = a_1,
\end{equation}
where $a_0, a_1, p, q$ are integers such that $\Delta = p^{2} + 4q > 0.$ Further, $\alpha = \frac{p + \sqrt{p^{2}+4q}}{2}$ and $\beta = \frac{p - \sqrt{p^{2}+4q}}{2}$ are roots of the characteristic equation. The Binet formula for this sequences are given by 
$$
w_{k} = \frac{A\alpha^{k}-B\beta^{k}}{\alpha - \beta},
$$

where $A=a_1-a_0\beta$ and $B=a_1-a_0\alpha$.\\

\bigskip

\section{Properties Of The Horadam Sequence}

\bigskip

\bigskip

In this section, we explore certain elementary properties of the Horadam sequence. Some these properties represents weighted sum formulas of this sequence. 

\begin{theorem}\label{2.1}
The Horadam sequence defined in \eqref{1.1} satisfies the following identities. 

\begin{itemize}
\item[\text{(a)}] $q^2w_k+pw_{k+3}=(p^2+q)w_{k+2},$\\
\item[\text{(b)}] $q^2w_k+w_{k+4}=(p^2+2q)w_{k+2},$\\
\item[\text{(c)}] $\left(q^3+p^2q^2\right)w_k +pw_{k+5}=\left(p^4+3p^2q+q^2\right)w_{k+2},$\\
\item[\text{(d)}] $(p^2q^2+2q^3)w_k+w_{k+6}=\left(p^4+3q^2 +4p^2q\right)w_{k+2},$\\
\item[\text{(e)}] $\left(p^4q^2+3p^2q^3+q^4\right)w_k +pw_{k+7}=\left(p^4+5p^4q+6p^2q^2+q^3\right)w_{k+2},$\\
\item[\text{(f)}] $\displaystyle\sum_{i=1}^{k} p^{k-i}q w_i=w_{k+2}-p^kw_2,$\\
\item[\text{(g)}] $\displaystyle\sum_{i=1}^{k} pq^{k-i} w_{2i-1}=w_{2k}-q^kw_0,$\\
\item[\text{(h)}] $\displaystyle\sum_{i=1}^{k} pq^{k-i} w_{2i}=w_{2k+1}-q^kw_1,$
\item[\text{(i)}] $w^{2}_{2k+3}-w^{2}_{2k+1}=w_{2k}w_{2k+6}-w_{2k}w_{2k+2}-
\frac{q^{2k} \cdot\left(qa^{2}_{0}+pa_{0}a_{1}-a^{2}_{1}\right)\cdot\left(p^3+pq+4q^3-p^2-4q\right)}{p^2+4q},$\\
\item[\text{(j)}] $w^{2}_{2k+3}+q^{2}w^{2}_{2k+1}=\left(p^{2}+2q\right)w_{2k}w_{2k+4}-
\frac{q^{2k} \cdot\left(qa^{2}_{0}+pa_{0}a_{1}-a^{2}_{1}\right)\cdot\left(p^3+p^2q^2+pq+8q^3\right)}{p^2+4q},$\\
\item[\text{(k)}] $w^{2}_{2k+3}-\left(p^{2}+2q\right)w^{2}_{2k+2}+q^{2}w^{2}_{2k+1}=\frac{q^{2k} \cdot\left(qa^{2}_{0}+pa_{0}a_{1}-a^{2}_{1}\right)\cdot\left(p^{6}+4p^{4}q+3p^{2}q^{2}-p^{3}-pq-8q^{3}\right)}{p^2+4q}.$
\end{itemize}
\end{theorem}
\begin{proof}
$\text{(a)}$ and $\text{(b)}$ can be proved using the recurrence relation of the sequence $w_k.$\\
We are proving $\text{(c)}.$ From $\text{(a)}$ we have $q^3w_k+pqw_{k+3}=(p^2q+q^2)w_{k+2}.$ From $\text{(b)}$ we have 
$p^2q^2w_k+p^2w_{k+4}=(p^4+2p^2q)w_{k+2}.$\\
Adding member with member the last two equalities, we obtain: 
$$\left(q^3+p^2q^2\right)w_k +p\left(pw_{k+4}+qw_{k+3}\right)=\left(p^4+3p^2q+q^2\right)w_{k+2}.$$ 
This last equality is equivalent to 
$$\left(q^3+p^2q^2\right)w_k +pw_{k+5}=\left(p^4+3p^2q+q^2\right)w_{k+2}.$$
We are proving $\text{(d)}.$ Using $\text{(b)}$ for $k\rightarrow k+2$ and then normally, we have:
$$q^2w_k+w_{k+6}=(p^2+2q)w_{k+4}-q^2w_{k+2}+q^2w_k=$$
$$=(p^2+2q)\cdot \left[ (p^2+2q)w_{k+2}-q^2w_{k}\right]-q^2w_{k+2}+q^2w_k=$$
$$=\left[(p^2+2q)^{2}-q^2\right]w_{k+2}+q^2(1-p^2-2q)w_k.$$
This implies that $(p^2q^2+2q^3)w_k+w_{k+6}=\left(p^4+3q^2 +4p^2q\right)w_{k+2}.$\\
The proof of $\text{(e)}$ is similar to the proof of $\text{(c)}$ and hence, it is omitted.\\
The proof of $\text{(f)}$ is based on mathematical induction.\\ 
Observe that the assertion in $\text{(f)}$ is true for $k=1$. Asuume that it is true for $k=t.$ Since, 
\begin{align*}
\sum_{i=1}^{t+1} p^{t+1-i}q w_i=&\sum_{i=1}^{t} p^{t+1-i}q w_i + q w_{t+1}\\
=& p\bigg(\sum_{i=1}^{t} p^{t-i}q w_i\bigg) + q w_{t+1}\\
=& p(w_{t+2}-p^{t}w_2)+ q w_{t+1}\\
=& w_{t+3}-p^{t+1}w_2,
\end{align*} 
 the assertion is true for $k=t+1.$ This proves $\text{(f)}$. The proofs of $\text{(g)}$ and $\text{(h)}$ are similar to that of $\text{(f)}$ and hence, are omitted.\\
$\text{(i)}$ Let $k,l$ be two positive integers, $k<l$. Using Binet formula for Horadam sequence, we have 
$$w_{2k}\cdot w_{2l}= \frac{A\alpha^{2k}-B\beta^{2k}}{\alpha - \beta}\cdot \frac{A\alpha^{2l}-B\beta^{2l}}{\alpha - \beta}=$$
$$=\frac{\left(A\alpha^{k+l}-B\beta^{k+l}\right)^{2}}{\left(\alpha - \beta\right)^{2}}+\frac{2AB\alpha^{k+l}\beta^{k+l}-AB\alpha^{2k}\beta^{2l}-AB\alpha^{2l}\beta^{2k} }{\left(\alpha - \beta\right)^{2}}=$$
$$=w^{2}_{k+l}+\frac{2AB\left(-q\right)^{k+l}-AB\left(-q\right)^{2k}\left(\alpha^{2l-2k}+\beta^{2l-2k}\right)}{\left(\alpha - \beta\right)^{2}}.$$
So, we obtained: 
\begin{equation*}
w_{2k}\cdot w_{2l}=w^{2}_{k+l}+\frac{2AB\left(-q\right)^{k+l}-ABq^{2k}\left(\alpha^{2l-2k}+\beta^{2l-2k}\right)}{\left(\alpha - \beta\right)^{2}} \tag{2.1.}
\end{equation*}
Making $l=k+1$ in the equality (2.1) and using Vi\`ete's relations for $\alpha$ and $\beta,$ we have:
$$w_{2k}\cdot w_{2k+2}=w^{2}_{2k+1}+\frac{-2ABq^{2k+1}-ABq^{2k}\left(\alpha^{2}+\beta^{2}\right)}{\left(\alpha - \beta\right)^{2}}= $$
$$=w^{2}_{2k+1}-\frac{ABq^{2k}\cdot\left(2q+p^{2}+2q\right)}{p^{2}+4q}=$$
$$=w^{2}_{2k+1}+q^{2k}\cdot\left(qa^{2}_{0}+pa_{0}a_{1}-a^{2}_{1}\right).$$
So, we obtained: 
\begin{equation*}
w_{2k}\cdot w_{2k+2}=w^{2}_{2k+1}+q^{2k}\cdot\left(qa^{2}_{0}+pa_{0}a_{1}-a^{2}_{1}\right) \tag{2.2.}
\end{equation*}
Making $l=k+2,$ respectively $l=k+3$ in the equality (2.1), using Vi\`ete's relations for $\alpha$ and $\beta$ and doing calculations similar to previous ones, we have:
\begin{equation*}
w_{2k}\cdot w_{2k+4}=w^{2}_{2k+2}+\frac{p^{2}\cdot q^{2k}\cdot\left(p^{2}+2q\right)\cdot\left(qa^{2}_{0}-a^{2}_{1}+pa_{0}a_{1}\right)}{p^{2}+4q} \tag{2.3.}
\end{equation*}
and
\begin{equation*}
w_{2k}\cdot w_{2k+6}=w^{2}_{2k+3}+\frac{q^{2k}\cdot\left(p^{3}+pq+4q^{3}\right)\cdot\left(qa^{2}_{0}-a^{2}_{1}+pa_{0}a_{1}\right)}{p^{2}+4q} \tag{2.4.}
\end{equation*}
From relations (2.2) and (2.4) we obtain:
$$w^{2}_{2k+3}-w^{2}_{2k+1}=w_{2k}w_{2k+6}-w_{2k}w_{2k+2}-
\frac{q^{2k} \cdot\left(qa^{2}_{0}+pa_{0}a_{1}-a^{2}_{1}\right)\cdot\left(p^3+pq+4q^3-p^2-4q\right)}{p^2+4q}.$$
$\text{(j)}$ From relations (2.2) and (2.4) we obtain:
$$w^{2}_{2k+3}+q^{2}w^{2}_{2k+1}=w_{2k}\cdot w_{2k+6}+q^{2}w_{2k}\cdot w_{2k+2}-\frac{q^{2k} \cdot\left(qa^{2}_{0}+pa_{0}a_{1}-a^{2}_{1}\right)\cdot\left(p^3+pq+4q^3+p^2q^2+4q^3\right)}{p^2+4q}.$$
Applying $\text{(b)}$ the last equality is equivalent with:
\begin{equation*}
w^{2}_{2k+3}+q^{2}w^{2}_{2k+1}=\left(p^{2}+2q\right)w_{2k}\cdot w_{2k+4}-\frac{q^{2k} \cdot\left(qa^{2}_{0}+pa_{0}a_{1}-a^{2}_{1}\right)\cdot\left(p^3+p^2q^2+pq+8q^3\right)}{p^2+4q} \tag{2.5.}
\end{equation*}
$\text{(k)}$ From relations (2.3) and (2.4) we obtain:
$$w^{2}_{2k+3}-\left(p^{2}+2q\right)w^{2}_{2k+2}+q^{2}w^{2}_{2k+1}=$$
$$=\frac{p^{2}\cdot q^{2k}\cdot\left(p^{2}+2q\right)^{2}\cdot\left(qa^{2}_{0}-a^{2}_{1}+pa_{0}a_{1}\right)}{p^{2}+4q} - \frac{q^{2k} \cdot\left(qa^{2}_{0}+pa_{0}a_{1}-a^{2}_{1}\right)\cdot\left(p^3+p^2q^2++pq+8q^3\right)}{p^2+4q}.$$
This last last equality is equivalent with:
$$w^{2}_{2k+3}-\left(p^{2}+2q\right)w^{2}_{2k+2}+q^{2}w^{2}_{2k+1}=\frac{q^{2k} \cdot\left(qa^{2}_{0}+pa_{0}a_{1}-a^{2}_{1}\right)\cdot\left(p^{6}+4p^{4}q+3p^{2}q^{2}-p^{3}-pq-8q^{3}\right)}{p^2+4q}.$$
\end{proof}

\begin{theorem}\label{2.2}
For every integer $n\geq0$ and for all $c\in \mathbb{R}-\{0\},$ 
$$c^{n+1}w_{n+1}=a_0+\sum_{i=0}^{n} c^{i}\bigg\{(p-1)w_{i}+(c-1)w_{i+1}+qw_{i-1}\bigg\}.$$
\end{theorem}
\begin{proof}
The proof is based on mathematical induction. For $n=0,$ both sides of the above equation are equal to $ca_1$. Assume that the assertion holds for $n=m$. Using the recurrence relation for the Horadam sequence, we get    
   \begin{align*}
   &a_0+\sum_{i=0}^{m+1} c^{i}\bigg\{(p-1)w_{i}+(c-1)w_{i+1}+qw_{i-1}\bigg\}\\
   &=a_0+\sum_{i=0}^{m} c^{i}\bigg\{(p-1)w_{i}+(c-1)w_{i+1}+qw_{i-1}\bigg\}\\
   &+c^{m+1}\bigg\{(p-1)w_{m+1}+(c-1)w_{m+2}+qw_{m}\bigg\}\\
   &=c^{m+1}w_{m+1}+c^{m+1}\bigg\{(p-1)w_{m+1}+(c-1)w_{m+2}+qw_{m}\bigg\}=c^{n+2}w_{m+2}.
   \end{align*}
   Thus, the assertion is true for $n=m+1.$ This completes the proof.
\end{proof}

As $(p,q)-$Fibonacci and $(p,q)-$Lucas sequence are special cases of Horadam sequence with initial values $a_0=0,a_1=1$ and $a_0=2,a_1=p$ respectively, they satisfy the identities in Theorem $2.1$.\\

The Binet formula for the $(p,q)-$Fibonacci and $(p,q)-$Lucas numbers are given by:
$$ \mathcal{F}_{p,q,k} = \frac{\alpha^{k}-\beta^{k}}{\alpha - \beta}~,~\mathcal{L}_{p,q,k} = \alpha^{k}+\beta^{k}.$$\\

 Using these Binet formulas, in the following theorem, we explore some properties of the $(p,q)-$Fibonacci and $(p,q)-$Lucas sequence.

\begin{theorem}\label{2.3}
Let $ \mathcal{F}_{p,q,k}$ and $ \mathcal{L}_{p,q,k}$ denote the $k^{th}$ term of the $(p,q)-$Fibonacci and $(p,q)-$Lucas sequence respectively. Then, we have the following identities.
\begin{itemize}
\item[\text{(a)}] $q\mathcal{F}_{p,q,k-1}+\mathcal{F}_{p,q,k+1}=\mathcal{L}_{p,q,k},$\\
\item[\text{(b)}] $q\mathcal{L}_{p,q,k-1}+\mathcal{L}_{p,q,k+1}=(p^2+4q)\mathcal{F}_{p,q,k},$\\
\item[\text{(c)}] $\mathcal{F}_{p,q,k+2}-q^2\mathcal{F}_{p,q,k-2}=p\mathcal{L}_{p,q,k},$\\
\item[\text{(d)}] $\mathcal{L}_{p,q,k+2}-q^2\mathcal{L}_{p,q,k-2}=p(p^2+4q)\mathcal{F}_{p,q,k},$\\
\item[\text{(e)}] $p\mathcal{F}_{p,q,k}+\mathcal{L}_{p,q,k}=2\mathcal{F}_{p,q,k+1},$\\
\item[\text{(f)}] $p\mathcal{L}_{p,q,k}+(p^2+4q)\mathcal{F}_{p,q,k}=2\mathcal{L}_{p,q,k+1},$\\
\item[\text{(g)}] $q^2\mathcal{F}_{p,q,k}+p\mathcal{L}_{p,q,k+2}=\mathcal{F}_{p,q,k+4},$\\
\item[\text{(h)}] $q^2\mathcal{L}_{p,q,k}+p(p^2+4q)\mathcal{F}_{p,q,k+2}=\mathcal{L}_{p,q,k+4},$\\
\item[\text{(i)}] $p\mathcal{F}_{p,q,k+2}+q\mathcal{L}_{p,q,k}=(p^2+2q)\mathcal{F}_{p,q,k+1},$\\
\item[\text{(j)}] $p\mathcal{L}_{p,q,k+2}+q(p^2+4q)\mathcal{F}_{p,q,k}=(p^2+2q)\mathcal{L}_{p,q,k+1},$\\
\item[\text{(k)}] $q^3\mathcal{F}_{p,q,k}+\mathcal{F}_{p,q,k+6}=(p^2+q)\mathcal{L}_{p,q,k+3},$\\
\item[\text{(l)}] $q^3\mathcal{L}_{p,q,k}+\mathcal{L}_{p,q,k+6}=(p^2+q)(p^2+4q)\mathcal{F}_{p,q,k+3},$\\
\item[\text{(m)}] $q^4\mathcal{F}_{p,q,k}+p\mathcal{F}_{p,q,k+8}=[(p^2+q)^2+pq(1+p+q)]\mathcal{F}_{p,q,k+4}.$\\
\end{itemize}
\end{theorem}
\begin{proof}
We prove $\text{(a)}$ only. Proofs of $\text{(b)}$-$\text{(m)}$ are similar. Using the Binet formula and the fact \\ $\alpha\beta=-q,$ we get
\begin{align*}
q\mathcal{F}_{p,q,k-1}+\mathcal{F}_{p,q,k+1}&=q\frac{\alpha^{k-1}-\beta^{k-1}}{\alpha - \beta}+\frac{\alpha^{k+1}-\beta^{k+1}}{\alpha - \beta}\\
&=\frac{1}{\alpha - \beta}\bigg[\alpha^{k}(\frac{q}{\alpha}+\alpha)-\beta^{k}(\frac{q}{\beta}+\beta)\bigg]\\
&=\frac{1}{\alpha - \beta}\bigg[\alpha^{k}(-\beta+\alpha)-\beta^{k}(-\alpha+\beta)\bigg]\\
&=\mathcal{L}_{p,q,k}.
\end{align*}
\end{proof}

The following theorem deals with the generating function and exponential generating function of a subsequence of the Horadam sequence. The proofs are similar to that of the Horadam sequence and hence are omitted.

\begin{theorem}\label{2.4}
For fixed integers $k$ and $m$ with $k>m\geq 0$ and $n \in \mathbb{N},$ 
\begin{itemize}
\item[\text{(a)}] the generating function of the Horadam number $w_{ kn+m}$ is $$
\sum_{n=0}^{\infty}w_{ kn+m} s^{n} = \frac{w_{ m} - (-q)^{k}w_{ m-k}s}{1-(\alpha^{k}+\beta^{k})s+(-q)^{k}s^{2}},
$$
\item[\text{(b)}] the exponential generating function of the Horadam number $w_{ kn+m}$ is
$$\sum\limits_{n=0}^{\infty}\frac{w_{ kn+m}}{n!}s^{n} = \frac{A\alpha^{m}e^{\alpha^{k}s} - B\beta^{m}e^{\beta^{k}s}}{\alpha-\beta}.$$
\end{itemize} 
\end{theorem}

The following theorem provides two sum formulas of a subsequence of the Horadam sequence.
\begin{theorem}\label{2.5}
For any natural numbers $m$ and $k$ with $k > m \geq 0,$ we have
\begin{enumerate}
\item[\text{(a)}] \begin{align*}
	\sum_{r=0}^{n} w_{ mr+k} =  \frac{(-q)^{m}w_{ mn+k}-w_{ mn+m+k}-(-q)^{m}w_{ k-m}+w_{ k}}{1+(-q)^{m}-(\alpha^m + \beta^m)},~~~~~~~~~~~~~~~~~~
	\end{align*}
	
	\item[\text{(b)}] \begin{align*}
	&\sum_{r=0}^{n} (-1)^r w_{ mr+k}\\
	&=  \frac{(-1)^{n+1}q^{m}w_{ mn+k}-(-1)^{n+1}w_{ mn+m+k}-q^{m}w_{ k-m}+w_{ k}}{1+(-q)^{m}-(-1)^m(\alpha^m + \beta^m)}.
	\end{align*}

\end{enumerate}

\end{theorem}
\begin{proof} We prove $\text{(a)}$ only. The proof of $\text{(b)}$ is similar.
Using the Binet formula for the Horadam sequence, we get
\begin{align*}
&\sum_{r=0}^{n} w_{ mr+k}\\
&= \sum_{r=0}^{n}\frac{A\alpha^{mr+k}-B\beta^{mr+k}}{\alpha - \beta}\\
&= \frac{1}{\alpha - \beta} \left[ A\alpha^{k}\left( \frac{\alpha^{mn+m} - 1}{\alpha^{m}-1}\right) - B\beta^{k}\left( \frac{\beta^{mn+m} - 1}{\beta^{m}-1}\right)\right] \\
&= \frac{1}{\alpha - \beta} \biggl[ \frac{(-q)^m (A\alpha^{mn+k}-B\beta^{mn+k}) - (A\alpha^{mn+m+k}-B\beta^{mn+m+k})}{1+(-q)^{m}-(\alpha^{m}+\beta^{m})} \\
&\hspace{3.5cm}- \frac{(-q)^m (A\alpha^{k-m}-B\beta^{k-m}) - (A\alpha^{k}-B\beta^{k})}{1+(-q)^{m}-(\alpha^{m}+\beta^{m})}\biggr].
\end{align*}

\end{proof}
 
 In the following theorem we explore weighted sum formulas for a subsequence of the Horadam sequence.
\begin{theorem}\label{2.6}
For any natural numbers $m$ and $k$ with $k > m \geq 0,$ we have
\begin{itemize}
		
	\item[\text{(a)}] \begin{align*}
	\sum_{r=0}^{n} rw_{ mr+k} = &n\bigg[\frac{(1+(-q)^{m})w_{ mn+m+k}-w_{ mn+2m+k}-(-q)^{m}w_{ mn+k}}{1+(-q)^{m}-(\alpha^m + \beta^m)}\\
	&\hspace{0.9cm}-\frac{w_{ mn+m+k}+q^{2m}w_{ mn+k-m}-2(-q)^{m}w_{ mn+k}}{(1+(-q)^{m}-(\alpha^m + \beta^m))^2}\bigg],
	\end{align*}
	
	\item[\text{(b)}] \begin{align*}
	&\sum_{r=0}^{n} (-1)^{r-1}rw_{ mr+k}\\
	 = &(-1)^{n+1} \bigg[\frac{(n-1-2n(-q)^m)w_{ mn+m+k}+q^{2m}(n-1)w_{ mn+k-m}}{(1+(-q)^{m}-(\alpha^m + \beta^m))^2}\\
	&\hspace{2.5cm}-\frac{nw_{ mn+2m+k}+(nq^{2m}-2(-q)^{m}(n-1))w_{ mn+k}}{(1+(-q)^{m}-(\alpha^m + \beta^m))^2}\bigg].
	\end{align*}
	
\end{itemize}

\end{theorem}

\begin{proof}
	These identities can be proved by using the proof of Theorem \ref{2.5} and the identity
	$$\sum_{r=1}^{n}rx^{r-1}=\frac{nx^{n+1}-(n+1)x^n+1}{(x-1)^2}, |x|<1$$
	which is obtained by differentiating the sum formula of a geometrical progression.
\end{proof}

In the following two theorems, we obtain a binomial weighted sum formula of the Horadam sequence.
\begin{theorem}\label{2.7}
For the integer $m\geq0$,
$\sum_{k=0}^{m} \begin{pmatrix}
m\\
k
\end{pmatrix} p^{k}w_{ k} q^{m-k} =  w_{ 2m}.$ 
\end{theorem}
\begin{proof}
Using the Binet formula of $w_k$, we get
\begin{align*}
&\sum_{k=0}^{m} \begin{pmatrix}
m\\
k
\end{pmatrix} p^{k}w_{ k} q^{m-k}\\
=& 
\sum_{k=0}^{m} \begin{pmatrix}
m\\
k
\end{pmatrix} p^{k} \frac{A\alpha^{k}-B\beta^{k}}{\alpha - \beta} q^{m-k} \\
=& \frac{A}{\alpha - \beta}\sum_{k=0}^{m} \begin{pmatrix}
m\\
k
\end{pmatrix} p^{k} \alpha^{k}q^{m-k}-\frac{B}{\alpha - \beta}\sum_{k=0}^{m} \begin{pmatrix}
m\\
k
\end{pmatrix} p^{k}\beta^{k} q^{m-k}\\
=& \frac{A}{\alpha - \beta} (p\alpha+q)^m -  \frac{B}{\alpha - \beta}(p\beta+q)^m
= \frac{A\alpha^{2m}-B\beta^{2m}}{\alpha - \beta}. 
\end{align*}
\end{proof}

\begin{theorem}\label{2.8}
Let $m$ be a non-negative integer. Then,
\begin{align*}
&\text{(a)}~ \sum_{n=0}^{m} \begin{pmatrix}
m\\
n
\end{pmatrix} w_{ 2n+k} q^{m-n} =  \left\{ \begin{array}{ll}
         w_{ k+m}\Delta^{\frac{m}{2}}, & \mbox{if $m ~even$};\\
        (A\alpha^{k+m}+B\beta^{k+m})\Delta^{\frac{m-1}{2}}, & \mbox{if $m ~odd$}.\end{array} \right. \\\\
&\text{(b)}~ \sum_{n=0}^{m} (-1)^{n} \begin{pmatrix}
m\\
n
\end{pmatrix} w_{ 2n+k} q^{m-n} =  \left\{ \begin{array}{ll}
         p^m w_{ k+m}, & \mbox{if $m ~even$};\\
        -p^m w_{ k+m}, & \mbox{if $m ~odd$}.\end{array} \right.\\\\
\end{align*}
\end{theorem}
\begin{proof}
The proof of this theorem is similar to that of Theorem \ref{2.7}. 
\end{proof}

\bigskip

\section{Properties Of The Horadam Symbol Elements}

\bigskip

\bigskip

In this section, we define Horadam symbol element and establish some properties for these elements. We also find the Binet formula and the generating function for these elements.\\

Consider the symbol algebra $S=\big(\frac{a,b}{K,\omega}\big)$ of degree $N$ generated by $x$ and $y$ over the field $K$. So $S$ has a $K$-basis\\

 $\{1,x,x^2,\cdots,x^{N-1},y,xy,x^2y,\cdots,x^{N-1}y,y^2,xy^2,x^2y^2,\cdots,x^{N-1}y^2,\cdots,$\\
 
 \hspace{8.2cm}$y^{N-1},xy^{N-1},\cdots,x^{N-1}y^{N-1}\}.$\\ 
  
We denote the $i^{th}$ term of the above basis as $e_{i-1}$ so that $\{e_0,e_1,\cdots,e_{N^2-1}\}$ is an ordered basis for S. We define the $k$-th term of Horadam symbol element as 
$$
W_{k} = \sum_{l=0}^{N^2-1} w_{k+l}e_{l},
$$ where $w_k$ is the $k^{th}$ Horadam number while $W_k$ is the $k^{th}$ Horadam symbol element. It is easy to see that the sequence of Horadam symbol elements satisfy the same recurrence
as Horadam sequence, that is, $$W_{k+1} = pW_{k} + qW_{k-1},$$ where the initial terms are $W_{0}=(w_0,w_1,\cdots,w_{n^2-1})$ and $W_{1}=(w_1,w_2,\cdots,w_{n^2})$. The auxiliary equation of this recurrence relation is given by $\lambda^2-p\lambda-q=0$ whose roots are $\alpha = \frac{p + \sqrt{p^{2}+4q}}{2}$ and $\beta = \frac{p - \sqrt{p^{2}+4q}}{2}.$ Let  $$
\underline{\alpha} = \sum_{l=0}^{N^2-1} \alpha^le_{l}
~\text{and}~
\underline{\beta} = \sum_{l=0}^{N^2-1} \beta^le_{l}.
$$

In the following theorem we obtain the Binet formula for the Horadam symbol elements.

\begin{theorem}\label{3.1}
The Binet formula of Horadam symbol elements is
$$W_k = \frac{A\underline{{\alpha}}\alpha^k - B\underline{{\beta}}\beta^k}{\alpha-\beta},$$ 
where $A = a_1 - a_0\beta, B = a_1 - a_0\alpha.$
\end{theorem}
\begin{proof}
Using the Binet formula for Horadam sequence, we get
\begin{align*}
W_{k} = \sum_{l=0}^{N^2-1} w_{k+l}e_{l}=& \sum_{l=0}^{N^2-1} \frac{A\alpha^{k+l} - B\beta^{k+l} }{\alpha-\beta} e_{l}\\
=& \frac{A\alpha^{k}}{\alpha-\beta} \sum_{l=0}^{N^2-1} \alpha^le_{l} - \frac{\beta^{k}}{\alpha-\beta} \sum_{l=0}^{N^2-1} \beta^le_{l}.
\end{align*}
\end{proof}

 In the following theorem, we prove two important identities for the Horadam symbol elements.
\begin{theorem}\label{3.2}
Let $W_k$ be the $k^{th}$ term of the Horadam symbol elements. Then the following identities hold.
\begin{itemize}
\item[\text{(a)}] $q^2W_k+pW_{k+3}=(p^2+q)W_{k+2},$\\
\item[\text{(b)}] $q^2W_k+W_{k+4}=(p^2+2q)W_{k+2},$\\
\item[\text{(c)}] $\left(q^3+p^2q^2\right)W_k +pW_{k+5}=\left(p^4+3p^2q+q^2\right)W_{k+2},$\\
\item[\text{(d)}] $(p^2q^2+2q^3)W_k+W_{k+6}=\left(p^4+3q^2 +4p^2q\right)W_{k+2},$\\
\item[\text{(e)}] $\left(p^4q^2+3p^2q^3+q^4\right)W_k +pW_{k+7}=\left(p^4+5p^4q+6p^2q^2+q^3\right)W_{k+2}.$
\end{itemize}
\end{theorem}
\begin{proof}
Using Theorem \ref{2.1} and the definition of Horadam symbol elements, we find
\begin{align*}
q^2W_k+pW_{k+3}=&q^2\sum_{l=0}^{N^2-1} w_{k+l}e_{l}+p\sum_{l=0}^{N^2-1} w_{(k+3)+l}e_{l}\\
=& \sum_{l=0}^{N^2-1}(q^2 w_{k+l}+pw_{(k+3)+l})e_{l}\\
=& \sum_{l=0}^{N^2-1} (p^2+q)w_{(k+2)+l}e_{l}.
\end{align*}
Thus $\text{(a)}$ is proved. The proofs of $\text{(b)}, \text{(c)}, \text{(d)}, \text{(e)}$ are similar to the proof of $\text{(a)}$ and hence, are omitted.

\end{proof}

In the next two theorems we establish some important sum formulas for the Horadam symbol elements.
\begin{theorem}\label{3.3}
The Horadam symbol elements satisfy
\begin{multicols}{2}
\begin{itemize}
\item[\text{(a)}] $\displaystyle\sum_{i=1}^{k} p^{k-i}q W_i=W_{k+2}-p^kW_2,$\\
\item[\text{(b)}] $\displaystyle\sum_{i=1}^{k} pq^{k-i} W_{2i-1}=W_{2k}-q^kW_0,$\\
\item[\text{(c)}] $\displaystyle\sum_{i=1}^{k} pq^{k-i} W_{2i}=W_{2k+1}-q^kW_1.$\\
\end{itemize}
\end{multicols}
\end{theorem}
\begin{proof}
We prove $\text{(a)}$ only. The proofs of $\text{(b)}$ and $\text{(c)}$ are similar. Using Theorem \ref{2.1} and the definition of Horadam symbol elements, we get
\begin{align*}
\sum_{i=1}^{k} p^{k-i}q W_i =& \sum_{i=1}^{k} p^{k-i}q\sum_{l=0}^{N^2-1} w_{i+l}e_{l}=\sum_{l=0}^{N^2-1} \bigg(\sum_{i=1}^{k} p^{k-i}q w_{i+l}\bigg)e_{l}\\
=& \sum_{l=0}^{N^2-1} \bigg(\sum_{i=1}^{k+l} p^{(k+l)-i}q w_{i}-p^k\sum_{i=1}^{l} p^{l-i}q w_{i}\bigg)e_{l}\\
=& \sum_{l=0}^{N^2-1} \bigg(w_{(k+2)+l}-p^kw_{2+l}\bigg)e_{l}\\
=& \bigg(\sum_{l=0}^{N^2-1} w_{(k+2)+l}e_{l}\bigg)-p^k\bigg(\sum_{l=0}^{N^2-1}w_{2+l}e_{l}\bigg).
\end{align*}
This proves $\text{(a)}$.
\end{proof}

\begin{theorem}\label{3.4}
For every integer $n\geq0$ and for all $c\in \mathbb{R}-\{0\},$  $$c^{n+1}W_{n+1}=W_0+\sum_{i=0}^{n} c^{i}\bigg\{(p-1)W_{i}+(c-1)W_{i+1}+qW_{i-1}\bigg\}.$$
\end{theorem}
\begin{proof}
The proof is similar to that of Theorem \ref{2.2} and hence, it is omitted.
\end{proof}

The two special cases of Horadam symbol elements are $(p,q)-$Fibonacci symbol elements $F_{p,q,k}$ and $(p,q)-$Lucas symbol elements $L_{p,q,k}$, defined by\\
 
$$F_{p,q,k}=\displaystyle\sum_{l=0}^{N^2-1} \mathcal{F}_{p,q,k+l}e_{l}~\text{and}~ L_{p,q,k}=\displaystyle\sum_{l=0}^{N^2-1} \mathcal{L}_{p,q,k+l}e_{l}$$ respectively, where $\mathcal{F}_{p,q,k} ~\text{and}~ \mathcal{L}_{p,q,k}$ denote the $(p,q)-$Fibonacci and $(p,q)-$Lucas number. The Binet formula for $(p,q)-$Fibonacci and $(p,q)-$Lucas symbol elements are respectively
$$F_{p,q,k} = \frac{\underline{{\alpha}}\alpha^k - \underline{{\beta}}\beta^k}{\alpha-\beta}~\text{and}~L_{p,q,k}=\underline{{\alpha}}\alpha^k + \underline{{\beta}}\beta^k.$$\\

Kilic\cite{Kilic} studied some binomial sums involving the $(p,q)-$Fibonacci numbers. The following theorem provides a similar sum formula for the Horadam symbol elements.\\

\begin{theorem}\label{3.5}
\begin{equation*}
W_{mn} = \mathcal{F}_{p,q,d}^{-n}\sum_{j=0}^{n} \begin{pmatrix}
n \\
j
\end{pmatrix} (-q)^{m(n-j)} \mathcal{F}_{p,q,d-m}^{n-j}  \mathcal{F}_{p,q,m}^{j}  W_{dj}.
\end{equation*}
\end{theorem}
\begin{proof}
Using the identity $$w_{mn+l} = \mathcal{F}_{p,q,d}^{-n} \sum_{j=0}^{n} \begin{pmatrix}
	n \\
	j
\end{pmatrix}  (-q)^{m(n-j)} \mathcal{F}_{p,q,d-m}^{n-j}  \mathcal{F}_{p,q,m}^{j} w_{dj+l},$$ the theorem follows directly.
\end{proof}

The following theorem shows that the $(p,q)-$Fibonacci and $(p,q)-$Lucas symbol elements satisfy similar identities as that of $(p,q)-$Fibonacci and $(p,q)-$Lucas numbers.

\begin{theorem}\label{3.6}
The $(p,q)-$Fibonacci and $(p,q)-$Lucas symbol elements satisfy 
\begin{itemize}
\item[\text{(a)}] $qF_{p,q,k-1}+F_{p,q,k+1}=L_{p,q,k},$\\
\item[\text{(b)}] $qL_{p,q,k-1}+L_{p,q,k+1}=(p^2+4q)F_{p,q,k},$\\
\item[\text{(c)}] $F_{p,q,k+2}-q^2F_{p,q,k-2}=pL_{p,q,k},$\\
\item[\text{(d)}] $L_{p,q,k+2}-q^2L_{p,q,k-2}=p(p^2+4q)F_{p,q,k},$\\
\item[\text{(e)}] $pF_{p,q,k}+L_{p,q,k}=2F_{p,q,k+1},$\\
\item[\text{(f)}] $pL_{p,q,k}+(p^2+4q)F_{p,q,k}=2L_{p,q,k+1},$\\
\item[\text{(g)}] $q^2F_{p,q,k}+pL_{p,q,k+2}=F_{p,q,k+4},$\\
\item[\text{(h)}] $q^2L_{p,q,k}+p(p^2+4q)F_{p,q,k+2}=L_{p,q,k+4},$\\
\item[\text{(i)}] $pF_{p,q,k+2}+qL_{p,q,k}=(p^2+2q)F_{p,q,k+1},$\\
\item[\text{(j)}] $pL_{p,q,k+2}+q(p^2+4q)F_{p,q,k}=(p^2+2q)L_{p,q,k+1},$\\
\item[\text{(k)}] $q^3F_{p,q,k}+F_{p,q,k+6}=(p^2+q)L_{p,q,k+3},$\\
\item[\text{(l)}] $q^3L_{p,q,k}+L_{p,q,k+6}=(p^2+q)(p^2+4q)F_{p,q,k+3},$\\
\item[\text{(m)}] $q^4F_{p,q,k}+F_{p,q,k+8}=[(p^2+q)^2+pq(1+p+q)]F_{p,q,k+4}.$\\
\end{itemize}
\end{theorem}
\begin{proof}
We prove $\text{(a)}$ only. The proofs of $\text{(b)}$-$\text{(m)}$ are similar. Using Theorem \ref{2.3} and the definition of $(p,q)-$Fibonacci and $(p,q)-$Lucas symbol elements, we find
\begin{align*}
qF_{p,q,k-1}+F_{p,q,k+1}=& q\sum_{l=0}^{N^2-1} \mathcal{F}_{p,q,(k-1)+l}e_{l} + \sum_{l=0}^{N^2-1} \mathcal{F}_{p,q,(k+1)+l}e_{l}\\
=& \sum_{l=0}^{N^2-1} [q\mathcal{F}_{p,q,(k+l)-1} + \mathcal{F}_{p,q,(k+l)+1}]e_{l}\\
=& \sum_{l=0}^{N^2-1} \mathcal{L}_{p,q,k+l}e_{l}= L_{p,q,k}.
\end{align*}
This proves $\text{(a)}.$
\end{proof}

 Generating functions are useful in solving linear homogeneous recurrences with constant coefficients. In the following theorem, we establish the generating function for the Horadam symbol elements. 

\begin{theorem}\label{3.7}
For fixed integers $k$ and $m$ with $k>m\geq 0$ and $n \in \mathbb{N},$ the generating functions of the subsequence $W_{ kn+m}$ of the Horadam symbol elements is given by
$$
\sum_{n=0}^{\infty}W_{ kn+m} s^{n} = \frac{W_{ m} - (-q)^{k}W_{ m-k}s}{1-(\alpha^{k}+\beta^{k})s+(-q)^{k}s^{2}}.
$$
\end{theorem}
\begin{proof}
Using Theorem \ref{2.4} and the definition of Horadam symbol elements, we find
\begin{align*}
\sum\limits_{n=0}^{\infty}W_{ kn+m}s^{n} =& \sum\limits_{n=0}^{\infty}\bigg[\sum_{l=0}^{N^2-1} w_{(kn+m)+l}e_{l}\bigg] s^{n}\\
=& \sum_{l=0}^{N^2-1}\bigg[\sum\limits_{n=0}^{\infty} w_{kn+(m+l)}s^{n}\bigg] e_{l}\\
=& \sum_{l=0}^{N^2-1}\bigg[\frac{w_{ m+l} - (-q)^{k}w_{ (m+l)-k}s}{1-(\alpha^{k}+\beta^{k})s+(-q)^{k}s^{2}}\bigg] e_{l}\\
=& \Bigg[\frac{\bigg(\displaystyle\sum_{l=0}^{N^2-1}w_{ m+l}e_{l}\bigg) - (-q)^{k}\bigg(\displaystyle\sum_{l=0}^{N^2-1}w_{ (m-k)+l}e_{l}\bigg)s}{1-(\alpha^{k}+\beta^{k})s+(-q)^{k}s^{2}}\Bigg]. 
\end{align*}
 
\end{proof}
\begin{corollary}\label{3.8}
 The generating function for the Horadam symbol elements $W_{n}$ is
$$
\sum_{n=0}^{\infty}W_{ n} s^{n} = \frac{W_{ 0} + (W_{ 1}-pW_0)s}{1-ps-qs^{2}}.
$$
\end{corollary}

 Exponential generating functions are also used
 for solving both homogeneous and nonhomogeneous linear recurrences with constant coefficients. In the following theorem, we establish the exponential generating function for the subsequences $W_{ kn+m}$ of the Horadam symbol elements.\\

\begin{theorem}\label{3.9}
For $m, n\in\mathbb{N}$, the exponential generating functions of the sub sequence $W_{ kn+m}$ of the Horadam symbol elements is given by
$$\sum\limits_{n=0}^{\infty}\frac{W_{ kn+m}}{n!}s^{n} = \frac{A\underline{\alpha}\alpha^{m}e^{\alpha^{k}s} - B\underline{\beta}\beta^{m}e^{\beta^{k}s}}{\alpha-\beta}.$$
\end{theorem}
\begin{proof}
Using Theorem \ref{2.4} and the definition of Horadam symbol elements, we get
\begin{align*}
\sum\limits_{n=0}^{\infty}\frac{W_{ kn+m}}{n!}s^{n}&=\sum\limits_{n=0}^{\infty}\bigg[\sum_{l=0}^{N^2-1} w_{(kn+m)+l}e_{l}\bigg]\frac{s^{n} }{n!}\\
&=\sum_{l=0}^{N^2-1}\bigg[\sum\limits_{n=0}^{\infty} w_{kn+(m+l)}e_{l}\bigg]\frac{s^{n} }{n!}\\
&=\sum_{l=0}^{N^2-1}\bigg[ \frac{A\alpha^{m+l}e^{\alpha^{k}s} - B\beta^{m+l}e^{\beta^{k}s}}{\alpha-\beta}e_{l}\bigg]\frac{s^{n} }{n!}\\
&=\Bigg[ \frac{A\alpha^{m}\bigg(\displaystyle\sum_{l=0}^{N^2-1}\alpha^{l}e_{l}\bigg)e^{\alpha^{k}s} - B\beta^{m}\bigg(\displaystyle\sum_{l=0}^{N^2-1}\beta^{l}e_{l}\bigg)e^{\beta^{k}s}}{\alpha-\beta}\Bigg]\frac{s^{n} }{n!}.
\end{align*}
\end{proof}
\begin{corollary}\label{3.10}
 The exponential generating function for the Horadam symbol elements $W_{n}$ is
$$\sum\limits_{n=0}^{\infty}\frac{W_{ n}}{n!}s^{n} = \frac{A\underline{\alpha}e^{\alpha s} - B\underline{\beta}e^{\beta s}}{\alpha-\beta}.$$
\end{corollary}

 Catalan's identity, Cassini formula and d'Ocagne's identity are usually obtained for the recurrent sequences where the multiplication is commutative. Since the elements of symbol algebra is associated with non-commutative multiplication, we will obtain two such identities in each case.
 
\begin{theorem}\label{3.11}
(Catlan Identity)Let $n, r \in \mathbb{Z}$, then
\begin{itemize}

\item[\text{(a)}] $$W_{ n-r}W_{ n+r}-W_{ n}^{2} = \frac{AB (-q)^{n-r}(\alpha^{r}-\beta^{r})(\beta^{r}\underline{\alpha}\underline{\beta}-\alpha^{r}\underline{\beta}\underline{\alpha})}{\Delta},$$

\item[\text{(b)}] $$W_{ n+r}W_{ n-r}-W_{ n}^{2} = \frac{AB (-q)^{n-r}(\alpha^{r}-\beta^{r})(\beta^{r}\underline{\beta}\underline{\alpha}-\alpha^{r}\underline{\alpha}\underline{\beta})}{\Delta}.$$
\end{itemize}

\end{theorem}
\begin{proof}
Using the Binet formula of Horadam symbol elements and the fact $\alpha\beta = -q$,  we get 
\begin{align*}
&W_{ n-r}W_{ n+r}-W_{ n}^{2}\\
&=  \left(\frac{A\underline{\alpha}\alpha^{n-r} - B\underline{\beta}\beta^{n-r}}{\alpha - \beta}\right) \left(\frac{A\underline{\alpha}\alpha^{n+r} - B\underline{\beta}\beta^{n+r}}{\alpha - \beta}\right)  - \left(\frac{A\underline{\alpha}\alpha^{n} - B\underline{\beta}\beta^{n}}{\alpha - \beta}\right)^2 \\
& = \frac{AB \underline{\alpha}\underline{\beta}(-q)^{n} \left( 1 - \frac{\beta^r}{\alpha^r} \right) + BA \underline{\beta}\underline{\alpha} (-q)^{n}  \left(1 - \frac{\alpha^r}{\beta^r} \right)}{(\alpha - \beta)^2}   \\
& = \frac{AB (-q)^{n-r}(\alpha^{r}-\beta^{r})(\beta^{r}\underline{\alpha}\underline{\beta}-\alpha^{r}\underline{\beta}\underline{\alpha})}{\Delta}.
\end{align*}
This proves $\text{(a)}$. In a similar way, $\text{(b)}$ can be proved.
\end{proof}
  
\begin{corollary}\label{3.12}
(Cassini Identity)Let $n \in \mathbb{Z}$, then
\begin{itemize}
\item[\text{(a)}] \begin{align*}
W_{ n-1}W_{ n+1}-W_{ n}^{2} = \frac{AB (-q)^{n-1}(\beta\underline{\alpha}\underline{\beta}-\alpha\underline{\beta}\underline{\alpha})}{\alpha-\beta},
\end{align*}

\item[\text{(b)}] \begin{align*}
W_{ n+1}W_{ n-1}-W_{ n}^{2} = \frac{AB (-q)^{n-1}(\beta\underline{\beta}\underline{\alpha}-\underline{\alpha}\underline{\beta})}{\alpha-\beta}.
\end{align*}

\end{itemize}
\end{corollary}

\begin{theorem} (d'Ocagne's identity).\label{3.13}
Let $m, n \in \mathbb N$ with $n \ge m$, then
\begin{itemize}
\item[\text{(a)}] $$W_{ n}W_{ m+1} - W_{ n+1}W_{ m} 
= \frac{(-q)^{m} AB (\alpha^{n-m}\underline{\alpha}\underline{\beta} - \beta^{n-m}\underline{\beta}\underline{\alpha})}{\alpha - \beta},$$

\item[\text{(b)}] $$W_{ m+1}W_{ n} - W_{ m}W_{ n+1} 
= \frac{(-q)^{m} AB (\alpha^{n-m}\underline{\beta}\underline{\alpha} - \beta^{n-m}\underline{\alpha}\underline{\beta})}{\alpha - \beta}.$$
\end{itemize}

\end{theorem}
\begin{proof}
Use of the Binet formula of Horadam symbol elements, yields

\begin{align*}
W_{ n}W_{ m+1} - W_{ n+1}W_{ m} &= \left(\frac{A\underline{\alpha}\alpha^{n} - B\underline{\beta}\beta^{n}}{\alpha - \beta}\right) \left(\frac{A\underline{\alpha}\alpha^{m+1} - B\underline{\beta}\beta^{m+1}}{\alpha - \beta}\right)\\
&- \left(\frac{A\underline{\alpha}\alpha^{n+1} - B\underline{\beta}\beta^{n+1}}{\alpha - \beta}\right) \left(\frac{A\underline{\alpha}\alpha^{m} - B\underline{\beta}\beta^{m}}{\alpha - \beta}\right) \\
&= \frac{ AB\underline{\beta}\underline{\alpha} \alpha^{m}\beta^{n}(\beta - \alpha) + AB \underline{\alpha}\underline{\beta} \alpha^{n}\beta^{m} (\alpha - \beta)}{(\alpha - \beta)^{2}} \\
&= \frac{(-q)^{m}AB\left( \alpha^{n-m}\underline{\alpha}\underline{\beta} - \beta^{n-m}\underline{\beta}\underline{\alpha} \right)}{\alpha - \beta}. 
\end{align*}
\end{proof}

Bolat and K\"{o}se\cite{Bolat} investigated certain properties of the $k-$Fibonacci sequence. The following theorem provide similar identities for the Horadam symbol elements.\\
 
\begin{theorem}\label{3.14}
	For arbitrary natural numbers $m$ and $n$, we get
	\begin{align*}
	W_{m}W_{n+1}+qW_{ m-1}W_{ n}=\frac{A^2\underline{{\alpha}}^2\alpha^{m+n} - B^2\underline{{\beta}}^2\beta^{m+n}}{\alpha-\beta}.
	\end{align*}
\end{theorem}
\begin{proof} Using the Binet formula for Horadam octonions, we have
	\begin{align*}
	W_{m}W_{n+1}+qW_{ m-1}W_{ n}&=\bigg(\frac{A\underline{{\alpha}}\alpha^{m} - B\underline{{\beta}}\beta^{m}}{\alpha-\beta}\bigg)\bigg( \frac{A\underline{{\alpha}}\alpha^{n+1}
	- B\underline{{\beta}}\beta^{n+1}}{\alpha-\beta}\bigg)\\
	&+ q\bigg(\frac{A\underline{{\alpha}}\alpha^{m-1} - B\underline{{\beta}}\beta^{m-1}}{\alpha-\beta}\bigg) \bigg( \frac{A\underline{{\alpha}}\alpha^{n} - B\underline{{\beta}}\beta^{n}}{\alpha-\beta}\bigg)\\
	&=\frac{A^2\underline{{\alpha}}^2\alpha^{m+n}(\alpha+\frac{q}{\alpha}) - B^2\underline{{\beta}}^2\beta^{m+n}(-\beta-\frac{q}{\beta})}{(\alpha-\beta)^2}\\
	&-\frac{AB\big(\underline{{\alpha}}\underline{{\beta}}\alpha^{m-1}\beta^{n}(\alpha\beta+q)+\underline{{\beta}}\underline{{\alpha}}\alpha^{n}\beta^{m-1}(\alpha\beta+q)\big)}{(\alpha-\beta)^2}\\
	&=\frac{A^2\underline{{\alpha}}^2\alpha^{m+n} - B^2\underline{{\beta}}^2\beta^{m+n}}{\alpha-\beta}.
	\end{align*}
	\end{proof}
	Observe that, \begin{align*}
	W_{n+1}W_{m}+qW_{ n}W_{ m-1}=\frac{A^2\underline{{\alpha}}^2\alpha^{m+n} - B^2\underline{{\beta}}^2\beta^{m+n}}{\alpha-\beta}.
	\end{align*}

In the following theorem, we establish two more such identities. 
\begin{theorem}\label{3.15}
	For any natural numbers $a,b,c,k$ and $a>k,b>k,c>k$, we have
	
	\begin{itemize}
	\item[\text{(a)}] 
	\begin{align*}
	\frac{W_{a}W_{b}-(-q)^kW_{a-k}W_{b-k}}{\mathcal{F}_{p,q,k}}=
	\frac{A^2\underline{{\alpha}}^2\alpha^{a+b-k} - B^2\underline{{\beta}}^2\beta^{a+b-k}}{\alpha-\beta}, 
	\end{align*}
	\item[\text{(b)}]
	\begin{align*}
		\frac{W_{a}W_{b}W_{c}-\mathcal{L}_{p,q,k}(-q)^kW_{a-k}W_{b-k}W_{c-k}+(-q)^{3k}W_{a-2k}W_{b-2k}W_{c-2k}}{\mathcal{L}_{p,q,k}\mathcal{F}_{p,q,k}^2}\hspace{0.5cm}\\
		=
		\frac{A^3\underline{{\alpha}}^3\alpha^{a+b+c-3k} - B^3\underline{{\beta}}^3\beta^{a+b+c-3k}}{\alpha-\beta}.
		\end{align*}
	\end{itemize} 	
where $\mathcal{F}_{p,q,k}$ and $\mathcal{L}_{p,q,k}$ denote the $k^{th}$ $(p,q)-$Fibonacci and $(p,q)-$Lucas numbers respectively.
\end{theorem}
\begin{proof} 
	The proof of $\text{(a)}$ and $\text{(b)}$ are similar to the proof of Theorem \ref{3.14} followed by the use of the identities 
	$$\alpha^{2k}-\mathcal{L}_{p,q,k}\alpha^{k}+(-q)^k=\beta^{2k}-\mathcal{L}_{p,q,k}\beta^{k}+(-q)^k=0$$ and
	 $$\alpha^{3k}-\mathcal{L}_{p,q,k}(-q)^k+\frac{(-q)^k}{\alpha^{3k}}=\beta^{3k}-\mathcal{L}_{p,q,k}(-q)^k+\frac{(-q)^k}{\beta^{3k}}=\mathcal{L}_{p,q,k}\mathcal{F}_{p,q,k}^2.$$ 
\end{proof}
In accordance with Theorem \ref{3.15}, we have the following two results.
\begin{theorem}\label{3.16}
For any natural numbers $a,b,c,d$ and $r$ with $a+b=c+d$, we have 
$$(-q)^r[W_{a}W_{b}-W_{c}W_{d}]=W_{a+r}W_{b+r}-W_{c+r}W_{d+r}$$	
\end{theorem}
\begin{theorem}\label{3.17}
	For any natural numbers $a,b,c,d,e,f$ and $r$ with $a+b+c=d+e+f$, we have
\begin{align*}
(-q)^{r}[\mathcal{L}_{p,q,k}(W_{a+r}W_{b+r}W_{c+r}-W_{d+r}W_{e+r}W_{f+r})-q^{2r}(W_{a}W_{b}W_{c}-W_{a}W_{b}W_{c})]\hspace{0.5cm}\\=W_{a+2r}W_{b+2r}W_{c+2r}-W_{d+2r}W_{e+2r}W_{f+2r}
\end{align*}
where $\mathcal{L}_{p,q,k}$ denote the $k^{th}$ $(p,q)-$Lucas number.
\end{theorem}
 
In the following theorem, we provide some sum and weighted sum formulas for a subsequence of the Horadam symbol elements.
 
\begin{theorem}\label{3.18}
For any natural numbers $m$ and $k$ with $k > m \geq 0,$ we have
\begin{itemize}
	\item[\text{(a)}] \begin{align*}
	\sum_{r=0}^{n} W_{ mr+k} =  \frac{(-q)^{m}W_{ mn+k}-W_{ mn+m+k}-(-q)^{m}W_{ k-m}+W_{ k}}{1+(-q)^{m}-(\alpha^m + \beta^m)},
	\end{align*}
	
	\item[\text{(b)}] \begin{align*}
	&\sum_{r=0}^{n} (-1)^r W_{ mr+k} \\
	&=  \frac{(-1)^{n+1}q^{m}W_{ mn+k}-(-1)^{n+1}W_{ mn+m+k}-q^{m}W_{ k-m}+W_{ k}}{1+(-q)^{m}-(-1)^m(\alpha^m + \beta^m)},
	\end{align*}
	
	\item[\text{(c)}] \begin{align*}
	&\sum_{r=0}^{n} rW_{ mr+k}\\
	 &= n\bigg[\frac{(1+(-q)^{m})W_{ mn+m+k}-W_{ mn+2m+k}-(-q)^{m}W_{ mn+k}}{1+(-q)^{m}-(\alpha^m + \beta^m)}\\
	&\hspace{3.0cm}-\frac{W_{ mn+m+k}+q^{2m}W_{ mn+k-m}-2(-q)^{m}W_{ mn+k}}{(1+(-q)^{m}-(\alpha^m + \beta^m))^2}\bigg],
	\end{align*}
	
	\item[\text{(d)}] \begin{align*}
	&\sum_{r=0}^{n} (-1)^{r-1}rW_{ mr+k}\\
	& = (-1)^{n+1}\bigg[\frac{(n-1-2n(-q)^m)W_{ mn+m+k}+q^{2m}(n-1)W_{ mn+k-m}}{(1+(-q)^{m}-(\alpha^m + \beta^m))^2}\\
	&\hspace{3.0cm}-\frac{nW_{ mn+2m+k}+(nq^{2m}-2(-q)^{m}(n-1))W_{ mn+k}}{(1+(-q)^{m}-(\alpha^m + \beta^m))^2}\bigg].
	\end{align*}
	
\end{itemize}

\end{theorem}
\begin{proof} Since all the proofs are similar, we prove $\text{(a)}$ only. 
Using the Theorem \ref{2.5}, Theorem \ref{2.6} and the definition of Horadam symbol elements, we get
\begin{align*}
&\sum_{r=0}^{n} W_{ mr+k}\\
&= \sum_{r=0}^{n}\bigg[\sum_{l=0}^{N^2-1} w_{(mr+k)+l}e_{l}\bigg]\\
&= \sum_{l=0}^{N^2-1}\bigg[\sum_{r=0}^{n} w_{mr+(k+l)}\bigg]e_{l}\\
&= \sum_{l=0}^{N^2-1}\bigg[\frac{(-q)^{m}w_{ mn+(k+l)}-w_{ mn+m+(k+l)}-(-q)^{m}w_{ (k+l)-m}+w_{ (k+l)}}{1+(-q)^{m}-(\alpha^m + \beta^m)} \bigg]e_{l}\\
&= \Bigg[\frac{(-q)^{m}\bigg(\displaystyle\sum_{l=0}^{N^2-1}w_{ (mn+k)+l}e_{l}\bigg)-\bigg(\displaystyle\sum_{l=0}^{N^2-1}w_{ (mn+m+k)+l}e_{l}\bigg)}{1+(-q)^{m}-(\alpha^m + \beta^m)} \\
&\hspace{3.5cm}- \frac{(-q)^{m}\bigg(\displaystyle\sum_{l=0}^{N^2-1}w_{ (k-m)+l}e_{l}\bigg)-\bigg(\displaystyle\sum_{l=0}^{N^2-1}w_{ k+l}e_{l}\bigg)}{1+(-q)^{m}-(\alpha^m + \beta^m)} \Bigg].
\end{align*}

\end{proof}
 
 The following corollary is an easy consequence of the above theorem.
\begin{corollary}\label{3.19} 
 For the Horadam symbol elements $W_n,$
 \begin{itemize}
	\item[\text{(a)}] \begin{align*}
	\sum_{r=0}^{n} W_{r}=\frac{(1-p)W_{0}+W_{1}+(-1)(qW_{n}+W_{n+1})}{1-(p+q)},
	\end{align*}
	
	\item[\text{(b)}] \begin{align*}
	\sum_{r=0}^{n} (-1)^r W_{ r}=\frac{(1+p)W_{ 0}-W_{1}+(-1)^{n+1}(qW_{n}-W_{n+1})}{1+(p-q)},
	\end{align*}
	
	\item[\text{(c)}] \begin{align*}
	\sum_{r=0}^{n} rW_{r}=n\bigg[\frac{(W_{n+1}-W_{n+2})+q(W_{n}-W_{n+1})}{1-(p+q)}-\frac{q^2W_{n-1}+2qW_{n}+W_{n+1}}{(1-(p^2+q))^2}\bigg],
	\end{align*}
	
	\item[\text{(d)}] \begin{align*}
	\sum_{r=0}^{n} (-1)^{r-1}rW_r=&(-1)^{n+1}\bigg[\frac{(n-1+2nq)W_{n+1}+q^{2}(n-1)W_{n-1}}{(1-(p+q))^2}\\
	&\hspace{2.8cm}-\frac{nW_{n+2}+(nq^{2}+2(q(n-1)))W_{n}}{(1-(p^2+q))^2}\bigg].~~
	\end{align*}
\end{itemize}
 \end{corollary}

 In the following theorem we provide some more identities involving Horadam symbol elements.
\begin{theorem}\label{3.20}
Let $m$ be a non-negative integer. Then,
\begin{align*}
&\text{(a)}~ \sum_{n=0}^{m} \begin{pmatrix}
m\\
n
\end{pmatrix} W_{ 2n+k} q^{m-n} =  \left\{ \begin{array}{ll}
         W_{ k+m}\Delta^{\frac{m}{2}}, & \mbox{if $m ~even$};\\
        (A\underline{\alpha}\alpha^{k+m}+B\underline{\beta}\beta^{k+m})\Delta^{\frac{m-1}{2}}, & \mbox{if $m ~odd$}.\end{array} \right. \\\\
&\text{(b)}~ \sum_{n=0}^{m} (-1)^{n} \begin{pmatrix}
m\\
n
\end{pmatrix} W_{ 2n+k} q^{m-n} =  \left\{ \begin{array}{ll}
         p^m W_{ k+m}, & \mbox{if $m ~even$};\\
        -p^m W_{ k+m}, & \mbox{if $m ~odd$}.\end{array} \right.\\\\
&\text{(c)}~ \sum_{n=0}^{m} \begin{pmatrix}
m\\
n
\end{pmatrix} W_{ n}W_{ n+k} q^{m-n} =  \left\{ \begin{array}{ll}
         (A^2\underline{\alpha}^{2}\alpha^{m+k} + B^{2}\underline{\beta}^{2}\beta^{m+k})\Delta^{\frac{m-2}{2}}, & \mbox{if $m ~even$};\\
        (A^2\underline{\alpha}^{2}\alpha^{m+k} - B^{2}\underline{\beta}^{2}\beta^{m+k})\Delta^{\frac{m-2}{2}}, & \mbox{if $m ~odd$}.\end{array} \right.\\\\ 
&\text{(d)}~ \sum_{n=0}^{m} \begin{pmatrix}
m\\
n
\end{pmatrix} p^{n}W_{ n} q^{m-n} =  W_{ 2m}.
\end{align*}
\end{theorem}
\begin{proof}
We prove $\text{(a)}$ only. Rest of the proofs are similar.\\

Using Binet formula of $W_{ k}$, we find
\begin{align*}
&\sum_{n=0}^{m} \begin{pmatrix}
m\\
n
\end{pmatrix} W_{ 2n+k} q^{m-n}\\
&= \sum\limits_{n=0}^{m} \begin{pmatrix}
m\\
n
\end{pmatrix} \left(\frac{A\underline{\alpha}\alpha^{2n+k} - B\underline{\beta}\beta^{2n+k}}{\alpha - \beta} \right) q^{m-n} \\
&= \frac{A\underline{\alpha}\alpha^{k}}{\alpha - \beta} \sum_{n=0}^{m} \begin{pmatrix}
m\\
n
\end{pmatrix} \left(\alpha^{2}\right)^{n} q^{m-n} - \frac{B\underline{\beta}\beta^{k}}{\alpha - \beta} \sum_{n=0}^{m} \begin{pmatrix}
m\\
n
\end{pmatrix} \left(\beta^{2}\right)^{n} q^{m-n}\\
&= \frac{A\underline{\alpha}\alpha^{k}(\alpha^{2}+q)^{m} - B\underline{\beta}\beta^{k}(\beta^{2}+q)^{m}}{\alpha - \beta} \\
&= \frac{A\underline{\alpha}\alpha^{k}(\alpha \sqrt{\Delta})^{m} - B\underline{\beta}\beta^{k}(- \beta\sqrt{\Delta})^{m}}{\alpha - \beta}\\
&= \bigg[\frac{A\underline{\alpha}\alpha^{k+m} +(-1)^{m+1} B\underline{\beta}\beta^{k+m}}{\alpha - \beta}\bigg]\Delta^{\frac{m}{2}}.
\end{align*}
\end{proof}

\smallskip

\section{Special Horadam symbol elements in symbol algebras of degree $3,$ over cyclotomic fields of finite fields}

\bigskip

Let $\{w_{k}(a_0,a_1,p,q) \},{k \geq 1}$ the Horadam sequence

\begin{equation}{\label{1.1}}
\hspace{3.0cm} w_{k+1} = pw_{k} + qw_{k-1};w_{0} = a_0,w_{1} = a_1,
\end{equation}
where $a_0, a_1, p, q$ are integers such that $\Delta = p^{2} + 4q > 0$ and $g.c.d\left(p, a_1,\right)$$>$$1.$
Let $r$ be a prime positive integer such that $r$ divides $g.c.d\left(p, a_1,\right),$ let the finite field $F=\mathbb{Z}_{r}$ and the field $K=\mathbb{Z}_{r}\left(\epsilon\right).$ Let $a,b$$\in$$K^{*}$ and $\epsilon$ be a primitive root of order $3$ of the unity. We consider the symbol algebra $S=\big(\frac{a,b}{K,\epsilon}\big)$ of degree $3$ generated by $x$ and $y$ over the field $K,$  with $x^3 = a$, $y^3 = b$, $yx = \epsilon xy.$
So $S$ has a $K$-basis $\left\{1,x,x^2,y,xy,x^{2}y,y^2,xy^{2},x^{2}y^{2}\right\}.$ 
We consider the reduced norm $\textbf{n}:$$S$$\longrightarrow$$K.$ If $z$$\in$$S$,
$z=\sum\limits_{i,j=0}^{2}x^{i}y^{j}c_{ij},$ the reduced norm of $z$ is
$$\textbf{n}\left(z\right)=a^{2}\cdot\left(c^{3}_{20}+ bc^{3}_{21}+ b^{2}c^{3}_{22}-3bc_{20}c_{21}c_{22}\right)
+ a\cdot\left(c^{3}_{10}+ bc^{3}_{11}+ b^{2}c^{3}_{12}-3bc_{10}c_{11}c_{12}\right)$$
$$-3a\cdot\left(c_{00}c_{10}c_{20} + bc_{01}c_{11}c_{21} + b^{2}c_{02}c_{12}c_{22}\right) - 
3ab\epsilon\cdot\left(c_{00}c_{12}c_{21} + c_{01}c_{10}c_{22} + c_{02}c_{11}c_{20}\right)-$$
\begin{equation}{\label{1.2}}
3ab\epsilon^{2}\cdot\left(c_{00}c_{11}c_{22} + c_{02}c_{10}c_{21} + c_{01}c_{12}c_{20}\right)+ c^{3}_{00}+ bc^{3}_{01}+ b^{2}c^{3}_{02}
-3bc_{00}c_{01}c_{02}.
\end{equation}
(See \cite[p.~299]{Pierce}).
The $k$-th Horadam symbol element in the symbol algebra $S$ is
$$W_{k} =w_{k}+ w_{k+1}x + w_{k+2}x^{2}+w_{k+3}y + w_{k+4}xy + w_{k+5}x^{2}y + w_{k+6}y^{2}+ w_{k+7}xy^{2}+ w_{k+8}x^{2}y^{2},$$ where $\left(w_k\right)_{k\geq0}$ is the Horadam sequence. In the following Proposition, we find the Horadam
symbol elements zero divisors in the symbol algebra $S$. For this, we use the properties of Horadam
symbol elements obtained in sections $2$ and $3.$
\begin{proposition}\label{4.1.} Let $W_{n}$ be the nth Horadam symbol element. Let $r$ be an odd prime positive integer, $r\neq 3$ such that $r$ $\mid$ $g.c.d\left(p, a_1\right),$ $r$ $\nmid$ $q,$ $r$ $\nmid$ $a_0,$ $\epsilon$ be
a primitive root of order $3$ of unity, let the field $K=\mathbb{Z}_{r}\left(\epsilon\right)$ and let $a,b$$\in$$\mathbb{Z}^{*},$ $r$ $\nmid$ $a,$ $r$ $\nmid$ $b.$ Then,\\
a) the norm of the symbol element $W_{2m}$ in the symbol algebra $S=\big(\frac{a,b}{K,\epsilon}\big)$ is:\\
$\textbf{n}\left(W_{2m}\right)=q^{3m}\cdot a^{3}_{0}\cdot \left[q^{3m}\cdot\left(a^{2}+a^{2}b^{2}+b^{2}+ab\right)+3abq^{2m}+1+
3abq^{2m}\epsilon\cdot \left(1-q^{m}\right)\right]$ (in $K$);\\
b) a Horadam symbol element $W_{2m}$  is a nontrivial zero divisor in the symbol algebra $S$ if and only if $q^{m}$$\equiv$$1$ (mod $r$) and $a^{2}+a^{2}b^{2}+b^{2}+4ab+1$$\equiv$$0$ (mod $r$).
\end{proposition}
\begin{proof} a)
 Since $r$ divides $p,$ applying Theorem 3.20 (d), we have $W_{2m}=W_{0}\cdot q^{m}$ (in $\mathbb{Z}_{r}$). Last equality is equivalent with 
$$W_{2m} = q^{m}\cdot \left(w_{0}+ w_{1}x + w_{2}x^{2}+w_{3}y + w_{4}xy + w_{5}x^{2}y + w_{6}y^{2}+ w_{7}xy^{2}+ w_{8}x^{2}y^{2}\right).$$ Applying Theorem 2.7, we have  $w_{2}=w_{4}=w_{6}=w_{8}=a_{0}\cdot q^{m}$ (in $\mathbb{Z}_{r}$). Using the recurrence relation of Horadam numbers, we have: $w_{3}=q\cdot a_{1}=\overline{0}$ $w_{5}=q^{2}\cdot a_{1}=\overline{0},$ $w_{7}=q^{3}\cdot a_{1}=\overline{0}$ (in $\mathbb{Z}_{r}$). So, it results that:
$$W_{2m} = q^{m}\cdot \left(a_{0}+ a_{0}q^{m}x^{2}+ a_{0}q^{m}xy+ a_{0}q^{m}y^{2}+ a_{0}q^{m}x^{2}y^{2}\right).$$
Using the relation (4.2), we obtain:
$$\textbf{n}\left(W_{2m}\right)=a^{2}\cdot a^{3}_{0}\cdot q^{6m}\cdot\left(1+b^{2}\right)
+ab\cdot a^{3}_{0}\cdot q^{6m}-3ab\epsilon\cdot a^{3}_{0}\cdot q^{6m}-3ab\epsilon^{2}\cdot a^{3}_{0}\cdot q^{5m}+ a^{3}_{0}\cdot q^{3m}+b^{2}\cdot a^{3}_{0}\cdot q^{6m}.$$
Using that $\epsilon^{2}=-1-\epsilon,$ the last equality becomes:\\
$\textbf{n}\left(W_{2m}\right)=q^{3m}\cdot a^{3}_{0}\cdot \left[q^{3m}\cdot\left(a^{2}+a^{2}b^{2}+b^{2}+ab\right)+3abq^{2m}+1+
3abq^{2m}\epsilon\cdot \left(1-q^{m}\right)\right]$ (in $K$).\\

b) First remark is that all Horadam symbol elements are $\neq 0.$\\
A Horadam symbol element $W_{2m}$ is a nontrivial zero divisor in the symbol algebra $S$ if and only if $\textbf{n}\left(W_{2m}\right)=0.$ This is equivalent with\\
$$q^{3m}\cdot a^{3}_{0}\cdot \left[q^{3m}\cdot\left(a^{2}+a^{2}b^{2}+b^{2}+ab\right)+3abq^{2m}+1+
3abq^{2m}\epsilon\cdot \left(1-q^{m}\right)\right]=0 \ (in \ K).\ $$
Since $r$ $\nmid$ $q,$ $r$ $\nmid$ $a_0,$ the last equality is equivalent with
$$q^{3m}\cdot\left(a^{2}+a^{2}b^{2}+b^{2}+ab\right)+3abq^{2m}+1+
3abq^{2m}\epsilon\cdot \left(1-q^{m}\right)=0 \ (in \  K).\ $$
This happens if and only if
$$\left\{
\begin{array}{c}
q^{3m}\cdot\left(a^{2}+a^{2}b^{2}+b^{2}+ab\right)+3abq^{2m}+1=\overline{0}\\
3abq^{2m}\cdot \left(1-q^{m}\right)=\overline{0}\ (in \  \mathbb{Z}_{r}).\
\end{array}%
\right.$$
Since $a,b$$\in$$\mathbb{Z}^{*},$  $r$ $\nmid$ $a,$ $r$ $\nmid$ $b,$ $r$ $\nmid$ $q$ and $\mathbb{Z}_{r}$ is a field, the second equation from the last system is equivalent with $q^{m}=\overline{1}$ (in $\mathbb{Z}_{r}$) $\Leftrightarrow$ $q^{m}$$\equiv$$1$ (mod $r$). \\
First equation of the system becomes $a^{2}+a^{2}b^{2}+b^{2}+4ab+1=\overline{0}$ (in $\mathbb{Z}_{r}$).
This is equivalent with $a^{2}+a^{2}b^{2}+b^{2}+4ab+1$$\equiv$$0$ (mod $r$).
\end{proof}

\begin{corollary}\label{4.2.} In the same hypotheses as in Proposition 4.1, if a Horadam symbol element $W_{2m}$ is a nontrivial zero divisor in the symbol algebra $S,$ then $r$$\equiv$$1$ (mod $4$).
\end{corollary}
\begin{proof} According to Proposition 4.1, if a Horadam symbol element $W_{2m}$ is a nontrivial zero divisor in the symbol algebra $S,$ then $a^{2}+a^{2}b^{2}+b^{2}+4ab+1$$\equiv$$0$ (mod $r$). This congruence is equivalent with $\left(a+b\right)^{2}+\left(ab+1\right)^{2}$$\equiv$$0$ (mod $r$). This implies that the Legendre symbol $\left(\frac{-\left(ab+1\right)^{2}}{r}\right)=1.$
So, the Legendre symbol $\left(\frac{-1}{r}\right)=1.$ This means $r$$\equiv$$1$ (mod $4$).
\end{proof}
\begin{proposition}\label{4.3.} Let $W_{n}$ be the nth Horadam symbol element. Let $r$ be an odd prime positive integer, $r\neq 3$ such that $r$ $\mid$ $g.c.d\left(p, a_1\right),$ $r$ $\nmid$ $q,$ $r$ $\nmid$ $a_0,$ $\epsilon$ be
a primitive root of order $3$ of unity, let the field $K=\mathbb{Z}_{r}\left(\epsilon\right)$ and let $a,b$$\in$$\mathbb{Z}^{*},$ $r$ $\nmid$ $a,$ $r$ $\nmid$ $b.$ Then, \\
a) the norm of the symbol element $W_{2m+1}$ in the symbol algebra $S=\big(\frac{a,b}{K,\epsilon}\big)$ is:\\
$$\textbf{n}\left(W_{2m+1}\right)=q^{3m+3}\cdot a^{3}_{0}\cdot \left(a^{2}bq^{6}+a+ ab^{2}q^{9}+ bq^{3}\right)\  (in\  \mathbb{Z}_{r}\subset K);$$
b) a Horadam symbol element $W_{2m+1}$  is a nontrivial zero divisor in the symbol algebra $S$ if and only if $a+ bq^{3}$$\equiv$$0$ (mod $r$) or $abq^{6}$$\equiv$$-1$ (mod $r$).
\end{proposition}
\begin{proof} a) The $2m+1$-th Horadam symbol element in the symbol algebra $S$ is
$$W_{2m+1} =w_{2m+1}+ w_{2m+2}x + w_{2m+3}x^{2}+w_{2m+4}y + w_{2m+5}xy + w_{2m+6}x^{2}y +$$
\begin{equation}{\label{1.3}}
+ w_{2m+7}y^{2}+ w_{2m+8}xy^{2}+ w_{2m+9}x^{2}y^{2}. 
\end{equation}
Applying Theorem 2.7, we have: $w_{2m+2}=a_{0}\cdot q^{m+1},$ $w_{2m+4}=a_{0}\cdot q^{m+2},$ $w_{2m+6}=a_{0}\cdot q^{m+3},$ $w_{2m+8}=a_{0}\cdot q^{m+4}$ (in $\mathbb{Z}_{r}$).

Using the fact that $r$ $\mid$ $g.c.d\left(p, a_1\right),$ and the recurrence relation of Horadam numbers, we have in $\mathbb{Z}_{r}$ : $w_{2m+1}=pw_{2m}+qw_{2m-1}=qw_{2m-1}=pqw_{2m-2}+q^{2}w_{2m-3}=q^{2}w_{2m-3}=...=q^{3}w_{2m-5}.$ Inductively, we obtain:
 $w_{2m+1}=q^{m}w_{1}=q^{m}a_{1}=\overline{0}.$ Analogously, we have: $w_{2m+3}=w_{2m+5}=w_{2m+7}=w_{2m+9}=\overline{0}$ (in $\mathbb{Z}_{r}$).\\
Using last relations, the equality (4.3) becomes:
\begin{equation}{\label{1.4}}
W_{2m+1} =a_{0}\cdot q^{m+1}\cdot \left(x+qy+q^{2}x^{2}y+q^{3}xy^{2}\right).
\end{equation}
From relations (4.2) and (4.4), we obtain:

$$\textbf{n}\left(W_{2m+1}\right)=q^{3m+9}\cdot a^{3}_{0}\cdot a^{2}\cdot b +a\cdot\left(q^{3m+3}\cdot a^{3}_{0}+ q^{3m+12}\cdot a^{3}_{0}\cdot b^{2}\right)+ q^{3m+6}\cdot a^{3}_{0}\cdot b \  (in\  \mathbb{Z}_{r}\subset K) \Leftrightarrow $$
\begin{equation}{\label{1.5}}
\textbf{n}\left(W_{2m+1}\right)=q^{3m+3}\cdot a^{3}_{0}\cdot \left(a^{2}bq^{6}+a+ ab^{2}q^{9}+ bq^{3}\right)\  (in\  \mathbb{Z}_{r}\subset K).
\end{equation}
b) A Horadam symbol element $W_{2m+1}$ is a nontrivial zero divisor in the symbol algebra $S$ if and only if $\textbf{n}\left(W_{2m+1}\right)=0.$ This is equivalent with\\
$$q^{3m+3}\cdot a^{3}_{0}\cdot \left(a^{2}bq^{6}+a+ ab^{2}q^{9}+ bq^{3}\right)=\overline{0}\  (in\  \mathbb{Z}_{r}\subset K).$$
Since $r$ $\nmid$ $q,$ $r$ $\nmid$ $a_0,$ the last equality becomes
$$a^{2}bq^{6}+a+ ab^{2}q^{9}+ bq^{3}=\overline{0}\  (in\  \mathbb{Z}_{r}).$$
The last equality is equivalent with
$$\left(a+ bq^{3}\right)\cdot \left(abq^{6}+1\right)=\overline{0}\  (in\  \mathbb{Z}_{r}).$$
Since $\mathbb{Z}_{r}$ is a field, we obtain that $a+ bq^{3}$$\equiv$$0$ (mod $r$) or $abq^{6}$$\equiv$$-1$ (mod $r$).
\end{proof}

\end{document}